\documentclass{article}%
\usepackage{amssymb}
\usepackage{graphicx}
\usepackage{amsmath}%
\usepackage{amsfonts}
\providecommand{\U}[1]{\protect\rule{.1in}{.1in}}
\providecommand{\U}[1]{\protect\rule{.1in}{.1in}}
\providecommand{\U}[1]{\protect\rule{.1in}{.1in}}
\providecommand{\U}[1]{\protect\rule{.1in}{.1in}}
\newtheorem{theorem}{Theorem}

\newtheorem{corollary}[theorem]{Corollary}

\newtheorem{lemma}[theorem]{Lemma}

\newtheorem{proposition}[theorem]{Proposition}
\newtheorem{remark}[theorem]{Remark}

\newenvironment{proof}[1][Proof]{\textbf{#1.} }{\ \rule{0.5em}{0.5em}}
\begin{document}

\title{The Jacobi Identity beyond Lie Algebras}
\author{Hirokazu Nishimura\\Institute of Mathematics, University of Tsukuba\\Tsukuba, Ibaraki, 305-8571\\Japan}
\maketitle

\begin{abstract}
Fr\"{o}licher and Nijenhuis recognized well in the middle of the previous
century that the Lie bracket and its Jacobi identity could and should exist
beyond Lie algebras. Nevertheless the conceptual status of their discovery has
been obscured by the genuinely algebraic techniques they exploited. The
principal objective in this paper is to show that the double dualization
functor in a cartesian closed category as well as synthetic differential
geometry provides an adequate framework, in which their discovery's conceptual
meaning appears lucid. The general Jacobi identity discovered by the author
[International Journal of Theoretical Physics, \textbf{36} (1997), 1099-1131]
will play a central role.

\end{abstract}

\section{Introduction}

Lie groups and their infinitesimal counterparts called Lie algebras were
introduced by Norwegian mathematician Sophus Lie in the 19th century. Lie
algebras are nonassociative algebras obeying the Jacobi identity instead. It
was Fr\"{o}licher and Nijenhuis (cf. \cite{fn} and \cite{n0}) in the middle of
the preceding century that realized the far-reaching nature of the Lie bracket
and its Jacobi identity (i.e., beyond Lie algebras) for the first time. They
have shown that tangent-vecor-valued differential forms enjoy a kind of Lie
bracket, which abides by a sort of the Jacobi identity. Nevertheless, because
of the genuinely algebraic techniques they used in order to establish their
marvelous discovery, the ubiquitous nature of the Lie bracket and its Jacobi
identity themselves has remained to be explored.

Now synthetic differential geometry, which is the avant-garde of differential
geometry, liberalizes ourselves. In particular, the general Jacobi identity
discovered by the author \cite{n1} more than a decade ago, which lies behind
the Jacobi identity of vector fields on a microlinear space, will play a
crucial role in this paper. For standard textbooks on synthetic differential
geometry the reader is referred to \cite{k2} or \cite{l1}.

The principal objective in this paper is to show that the double dualization
functor in a cartesian closed category as well as the general Jacobi identity
established in synthetic differential geometry provides us with the desired
framework. Our approach is completely combinatorial or geometric in sharp
contrast to Fr\"{o}licher and Nijenhuis' genuinely algebraic approach. After
some preliminaries, we present our discovery in the most abstract form in
\S \ref{3}. This abstract Jacobi identity for the double dualization functor
is then specialized in two distinct ways. In \S \ref{4} we specialize the
abstract Jacobi identity to tangent-vector-valued differential forms, while we
do so for Schwartz distributions in \S \ref{5}.

Last but not least, I gladly acknowledge my indebtedness to Professor Anders
Kock (Aarhus University), who kindly helped me pay due attention to the double
dualization functor in a cartesian closed category. His sincere and detailed
advice has improved the previous paper \cite{n4} considerably.

\section{Preliminaries\label{2}}

\subsection{The Double Dualization Functor\label{2.3}}

Let $\mathcal{E}$ be a cartesian closed category. It is well known that
cartesian closed categories and typed $\lambda$-calculi are essentially
equivalent, for which the reader is referred to, e.g., Chapter 4 of \cite{ama}
or Chapter 6 of \cite{bar}, so that we can speak about $\mathcal{E}$ in terms
of typed $\lambda$-calculi. Given two objects $A,B$ in $\mathcal{E} $, we
denote by $[A\rightarrow B]$ the exponential of $A$ over $B$, which is often
written $B^{A}$. We now fix an object $M$ in $\mathcal{E}$, which gives rise
to the \textit{double dualization functor} assigning $[[A\rightarrow
M]\rightarrow M]$ to each object \ $A$ in $\mathcal{E}$. Given $f\in
\lbrack\lbrack A\rightarrow M]\rightarrow M]$ and $g\in\lbrack\lbrack
B\rightarrow M]\rightarrow M]$, two kinds of \textit{convolution} of $f$ and
$g$, both of which belong to $[[A\times B\rightarrow M]\rightarrow M] $, are
defined, as is familiar in the theory of distributions, to be
\begin{align*}
\lambda h  &  \in[A\times B\rightarrow M].f(\lambda a\in A.g(\lambda b\in
B.h(a,b)))\\
\lambda h  &  \in[A\times B\rightarrow M].g(\lambda b\in B.f(\lambda a\in
A.h(a,b)))
\end{align*}
The former is denoted by $f\ast g$, while the latter is denoted by
$f\widetilde{\ast}g$. By identifying $A\times B$ and $B\times A$ naturally, we
can say that $f\widetilde{\ast}g$ is no other than $g\ast f$. It should be
obvious that

\begin{lemma}
\label{t2.3.1}Given $f\in\lbrack\lbrack A\rightarrow M]\rightarrow M]$,
$g\in\lbrack\lbrack B\rightarrow M]\rightarrow M]$ and $h\in\lbrack\lbrack
C\rightarrow M]\rightarrow M]$, we have
\begin{align*}
(f\ast g)\ast h  &  =f\ast(g\ast h)\\
(f\widetilde{\ast}g)\widetilde{\ast}h  &  =f\widetilde{\ast}(g\widetilde{\ast
}h)
\end{align*}

\end{lemma}

This lemma enables us to write, e.g., $f\ast g\ast h$ without parentheses in
place of $(f\ast g)\ast h$ or $f\ast(g\ast h)$.

If $a$ is an element of $A$ (i.e., $a$ is a global section $1\rightarrow A$),
then $\lambda f\in\lbrack A\rightarrow M].f(a)$ is denoted by $\delta_{a}$
and, exploiting the terminology in the theory of distributions, is called the
\textit{Dirac distribution at} $a$. The following lemma should be obvious.

\begin{lemma}
\label{t2.3.2}If one of $f\in\lbrack\lbrack A\rightarrow M]\rightarrow M]$ and
$g\in\lbrack\lbrack B\rightarrow M]\rightarrow M]$ is a Dirac distribution,
then $f\ast g$ and $f\widetilde{\ast}g$ coincide.
\end{lemma}

If both $A$ and $B$ are a terminal object $1$, then both $[[A\rightarrow
M]\rightarrow M]$ and $[[B\rightarrow M]\rightarrow M]$ can naturally be
identified with $[M\rightarrow M]$, so that the above two convolutions
degenerate into the composition of mappings in such a way that
\begin{align*}
f\ast g  &  =f\circ g\\
f\widetilde{\ast}g  &  =g\circ f
\end{align*}

\subsection{Synthetic Differential Geometry\label{2.0}}

We assume that the reader is familiar with Lavendhomme's textbook \cite{l1} on
synthetic differential geometry up to Chapter 4. From now on our discussion
will be done within an adequate universe of synthetic differential geometry,
as in Lavendhomme's textbook \cite{l1}. We denote by $D$ the subset of
$\mathbb{R}$ (the extended set of real numbers satisfying the generalized
Kock-Lawvere axiom so that $\mathbb{R}$ is microlinear) consisting of elements
$d$ of $\mathbb{R}$ with $d^{2}=0$. \ We shall let $M$ and $N$ with or without
subscripts denote microlinear spaces in the sense of Definition 1 in \S 2.3 of
\cite{l1}.

Given $\gamma\in\lbrack D^{p}\rightarrow M]$, $\alpha\in\mathbb{R}$ and a
natural number $i$ with $1\leq i\leq p$, we define $\alpha\underset{i}{\cdot
}\gamma\in\lbrack D^{p}\rightarrow M]$ to be
\[
(\alpha\underset{i}{\cdot}\gamma)(d_{1},...,d_{i-1},d_{i},d_{i+1}%
,...,d_{p})=\gamma(d_{1},...,d_{i-1},\alpha d_{i},d_{i+1},...,d_{p})
\]
for any $(d_{1},...,d_{i-1},d_{i},d_{i+1},...,d_{p})\in D^{p}$. We write
$\mathbb{S}_{p}$ for the permutation group of the first $p$ natural numbers,
namely, $1,...,p$. Given $\sigma\in\mathbb{S}_{p}$, we denote by
$\varepsilon_{\sigma}$ its signature. Given $\gamma\in\lbrack D^{p}\rightarrow
M]$ and $\sigma\in\mathbb{S}_{p}$, we define $\gamma^{\sigma}\in\lbrack
D^{p}\rightarrow M]$ to be
\[
\gamma^{\sigma}(d_{1},...,d_{p})=\gamma(d_{\sigma(1)},...,d_{\sigma(p)})
\]
for any $(d_{1},...,d_{p})\in D^{p}$. Given $\varphi\in\lbrack\lbrack
D^{p}\rightarrow M]\rightarrow M]$ and $\alpha\in\mathbb{R}$, we define
$\alpha\underset{i}{\cdot}\varphi\in\lbrack\lbrack D^{p}\rightarrow
M]\rightarrow M]$ ($1\leq i\leq p$) to be
\[
(\alpha\underset{i}{\cdot}\varphi)(\gamma)=\varphi(\alpha\underset{i}{\cdot
}\gamma)
\]
for any $\gamma\in M^{D^{p}}$. Given $\varphi\in\lbrack\lbrack D^{p}%
\rightarrow M]\rightarrow M]$ and any $\sigma\in\mathbb{S}_{p}$, we define
$\varphi^{\sigma}\in\lbrack\lbrack D^{p}\rightarrow M]\rightarrow M]$ to be
\[
\varphi^{\sigma}(\gamma)=\varphi(\gamma^{\sigma})
\]
for any $\gamma\in\lbrack D^{p}\rightarrow M]$. Given $\varphi\in
\lbrack\lbrack D^{p}\rightarrow M]\rightarrow M]$ and $\sigma,\tau
\in\mathbb{S}_{p}$, it is easy to see that
\[
\varphi^{\sigma\tau}(\gamma)=\varphi(\gamma^{\sigma\tau})=\varphi
((\gamma^{\sigma})^{\tau})=\varphi^{\tau}(\gamma^{\sigma})=(\varphi^{\tau
})^{\sigma}(\gamma)
\]
for any $\gamma\in\lbrack D^{p}\rightarrow M]$, so that $\varphi^{\sigma\tau
}=(\varphi^{\tau})^{\sigma}$.

\subsection{Vector Fields\label{2.1}}

In synthetic differential geometry vector fields on $M$ can be viewed in three
distinct but equivalent ways, which is based upon the following familiar
exponential laws:
\begin{align*}
\lbrack M  &  \rightarrow[D\rightarrow M]]\\
&  =[M\times D\rightarrow M]\\
&  =[D\rightarrow\lbrack M\rightarrow M]]
\end{align*}
The first viewpoint, which is based upon the first exponential form in the
above and is highly orthodox in traditional differential geometry, is to
regard a vector field on $M$ as a section of the canonical projection
$[D\rightarrow M]\rightarrow M$. The second viewpoint, which is based upon the
middle exponential form in the above, is to look upon a vector field on $M$ as
an infinitesimal flow on $M$. The third viewpoint, which is most radical and
is based upon the last exponential form in the above, is to speak of a vector
field on $M$ as an infinitesimal transformation of $M$. For the detailed
exposition of these three viewpoints on vector fields and their equivalence,
the reader is referred to \S 3.2 of \cite{l1}.

\subsection{The General Jacobi Identity\label{2.2}}

The notion of strong difference $\overset{\cdot}{-}$ was introduced by Kock
and Lavendhomme \cite{kl} into synthetic differential geometry. The notion of
strong difference $\overset{\cdot}{-}$ can be relativized. Since
$[D^{3}\rightarrow M]=[D^{2}\rightarrow\lbrack D\rightarrow M]]$, microcubes
on $M$ can be viewed as microsquares on $[D\rightarrow M]$. According to which
$D$ in the right-hand side of $D^{3}=D\times D\times D$ appears in the
subformula $[D\rightarrow M]$ of $[D^{2}\rightarrow\lbrack D\rightarrow M]]$,
we get the three relativized strong differences $\underset{i}{\overset{\cdot
}{-}}$ $(i=1,2,3)$, for which we have the following general Jacobi identity.

\begin{theorem}
Let $\gamma_{123},\gamma_{132},\gamma_{213},\gamma_{231},\gamma_{312}%
,\gamma_{321}\in\lbrack D^{3}\rightarrow M]$. As long as the following three
expressions are well defined, they sum up only to vanish:
\begin{align*}
&  (\gamma_{123}\overset{\cdot}{\underset{1}{-}}\gamma_{132})\overset{\cdot
}{-}(\gamma_{231}\overset{\cdot}{\underset{1}{-}}\gamma_{321})\\
&  (\gamma_{231}\overset{\cdot}{\underset{2}{-}}\gamma_{213})\overset{\cdot
}{-}(\gamma_{312}\overset{\cdot}{\underset{2}{-}}\gamma_{132})\\
&  (\gamma_{312}\overset{\cdot}{\underset{3}{-}}\gamma_{321})\overset{\cdot
}{-}(\gamma_{123}\overset{\cdot}{\underset{3}{-}}\gamma_{213})
\end{align*}

\end{theorem}

The theorem was established by the author in \cite{n1} and has been reproved
twice by himself in \cite{n2} and \cite{n3}, where K. Osoekawa aided the
author in computer algebra in the latter paper. The Jacobi identity of vector
fields on $M$ follows from the above theorem at once, as was noted in
\cite{n1}.

\section{The Jacobi Identity for the Double Dualization Functor\label{3}}

Let $n$ be a natural number. Let $A$ be a space with $a_{0}\in A$. An
$n$\textit{-dimensional }$(A,a_{0})$\textit{-icon on }$M$ or, more simply, an
$(A,a_{0})$\textit{-}$n$\textit{-icon on} $M$ is simply a mapping $\xi
:D^{n}\rightarrow\lbrack\lbrack A\rightarrow M]\rightarrow M]$ with
$\xi(0,...,0)=\delta_{a_{0}}$. Let $B$ be another object in $\mathcal{E}$ with
$b_{0}\in B$. Let $m$ be a natural number. Given an $(A,a_{0})$-$n$-icon
$\xi_{1}$on $M$ and a $(B,b_{0})$-$m$-icon $\xi_{2}$ on $M$, their
compositions $\xi_{1}\circledast\xi_{2}$ and $\xi_{1}\widetilde{\circledast
}\xi_{2}$, both of which are $(A\times B,(a_{0},b_{0}))$-$(m+n)$-icons, are
defined to be
\begin{align*}
(\xi_{1}\circledast\xi_{2})(\underline{d_{1}},\underline{d_{2}})  &  =\xi
_{1}(\underline{d_{1}})\ast\xi_{2}(\underline{d_{2}})\\
(\xi_{1}\widetilde{\circledast}\xi_{2})(\underline{d_{1}},\underline{d_{2}})
&  =\xi_{1}(\underline{d_{1}})\widetilde{\ast}\xi_{2}(\underline{d_{2}})
\end{align*}
for any $(\underline{d_{1}},\underline{d_{2}})\in D^{n}\times D^{m}=D^{m+n}$.
In particular, if $m=n=1$, then, by Lemma \ref{t2.3.2}, we have $\xi
_{1}\circledast\xi_{2}\mid_{D(2)}=\xi_{1}\widetilde{\circledast}\xi_{2}%
\mid_{D(2)}$, so that we can define their strong difference, called their
\textit{Lie bracket }and denoted by $\left\lfloor \xi_{1},\xi_{2}\right\rfloor
$, to be
\[
\left\lfloor \xi_{1},\xi_{2}\right\rfloor =\xi_{1}\widetilde{\circledast}%
\xi_{2}\overset{\cdot}{-}\xi_{1}\circledast\xi_{2}%
\]

Let $A$ and $B$ be objects in $\mathcal{E}$ with $a_{0}\in A$ and $b_{0}\in B
$.

We will show that the bracket $\left\lfloor ,\right\rfloor $ is antisymmetric.

\begin{theorem}
\label{t3.1}Let $A$ and $B$ be spaces with $a_{0}\in A$ and $b_{0}\in B$. Let
$\xi_{1}$ be an $(A,a_{0})$-$1$-icon on $M$ and $\xi_{2}$ a $(B,b_{0})$%
-$1$-icon on $M$. Then we have the following antisymmetry:
\[
\left\lfloor \xi_{1},\xi_{2}\right\rfloor +\left\lfloor \xi_{2},\xi
_{1}\right\rfloor =0
\]

\end{theorem}

\begin{remark}
Since $A\times B$ and $B\times A$ can naturally be identified, not only
$\left\lfloor \omega_{1},\omega_{2}\right\rfloor $ but also $\left\lfloor
\omega_{2},\omega_{1}\right\rfloor $ is to be regarded as a mapping
$D\rightarrow\lbrack\lbrack A\times B\rightarrow M]\rightarrow M]$. The reader
should be aware that the permutation sigma, which shall be omitted
intentionally for notational simplicity in this section, should be inserted in
the identification of $A\times B$ and $B\times A$ from a very strict
viewpoint. This comment should be recalled in the next section.
\end{remark}

\begin{proof}
This follows from Propositions 4 and 6 in \S 3.4 of Lavendhomme \cite{l1}.
More specifically we have
\begin{align*}
&  \left\lfloor \xi_{1},\xi_{2}\right\rfloor +\left\lfloor \xi_{2},\xi
_{1}\right\rfloor \\
&  =(\xi_{1}\widetilde{\circledast}\xi_{2}\overset{\cdot}{-}\xi_{1}%
\circledast\xi_{2})+(\xi_{2}\widetilde{\circledast}\xi_{1}\overset{\cdot}%
{-}\xi_{2}\circledast\xi_{1})\\
&  =(\xi_{1}\widetilde{\circledast}\xi_{2}\overset{\cdot}{-}\xi_{1}%
\circledast\xi_{2})+(\xi_{1}\circledast\xi_{2}\overset{\cdot}{-}\xi
_{1}\widetilde{\circledast}\xi_{2})\\
&  \text{[By Proposition 6 in \S 3.4 of Lavendhomme \cite{l1},}\\
&  \text{since }\xi_{2}\widetilde{\circledast}\xi_{1}\text{ can be identified
with }\\
(d_{1},d_{2})  &  \in D^{2}\mapsto\xi_{1}(d_{2})\ast\xi_{2}(d_{1})\in
\lbrack\lbrack A\times B,M],M]\text{, and}\\
&  \text{similarly for }\xi_{2}\circledast\xi_{1}\text{]}\\
&  =0\text{ \ }\\
&  \text{[By Proposition 4 in \S 3.4 of Lavendhomme \cite{l1}]}%
\end{align*}

\end{proof}

\begin{theorem}
\label{t3.2}Let $A$, $B$ and $C$ be objects in $\mathcal{E}$ with $a_{0}\in A
$, $b_{0}\in B$ and $c_{0}\in C$. Let $\xi_{1}$ be an $(A,a_{0})$-$1$-icon on
$M$, $\xi_{2}$ a $(B,b_{0})$-$1$-icon on $M$, and $\xi_{3}$ a $(C,c_{0}) $%
-$1$-icon on $M$. Then we have the following Jacobi identity:
\[
\left\lfloor \xi_{1},\left\lfloor \xi_{2},\xi_{3}\right\rfloor \right\rfloor
+\left\lfloor \xi_{2},\left\lfloor \xi_{3},\xi_{1}\right\rfloor \right\rfloor
+\left\lfloor \xi_{3},\left\lfloor \xi_{1},\xi_{2}\right\rfloor \right\rfloor
=0
\]

\end{theorem}

\begin{remark}
As in Theorem \ref{t3.1}, not only $\left\lfloor \xi_{1},\left\lfloor \xi
_{2},\xi_{3}\right\rfloor \right\rfloor $ but also both $\left\lfloor \xi
_{2},\left\lfloor \xi_{3},\xi_{1}\right\rfloor \right\rfloor $ and
$\left\lfloor \xi_{3},\left\lfloor \xi_{1},\xi_{2}\right\rfloor \right\rfloor
$ are to be regarded as mappings $D\rightarrow\lbrack\lbrack A\times B\times
C\rightarrow M]\rightarrow M]$.
\end{remark}

In order to establish this theorem, we need the following simple lemma, which
is a tiny generalization of Proposition 2.6 of \cite{n1}.

\begin{lemma}
\label{t3.3}Let $\xi$ be an $(A,a_{0})$-$1$-icon on $M$, and $\xi_{1}$ and
$\xi_{2}$ $(B,b_{0})$-$2$-icons on $M$ with $\xi_{1}\mid_{D(2)}=\xi_{2}%
\mid_{D(2)}$. Then the following formulas are both meaningful and valid.
\begin{align*}
\xi\circledast\xi_{1}\underset{1}{\overset{\cdot}{-}}\xi\circledast\xi_{2}  &
=\xi\circledast(\xi_{1}\overset{\cdot}{-}\xi_{2})\\
\xi\widetilde{\circledast}\xi_{1}\underset{1}{\overset{\cdot}{-}}\xi
\widetilde{\circledast}\xi_{2}  &  =\xi\widetilde{\circledast}(\xi_{1}%
\overset{\cdot}{-}\xi_{2})\\
\xi_{1}\circledast\xi\underset{3}{\overset{\cdot}{-}}\xi_{2}\circledast\xi &
=(\xi_{1}\overset{\cdot}{-}\xi_{2})\circledast\xi\\
\xi_{1}\widetilde{\circledast}\xi\underset{3}{\overset{\cdot}{-}}\xi
_{2}\widetilde{\circledast}\xi &  =(\xi_{1}\overset{\cdot}{-}\xi
_{2})\widetilde{\circledast}\xi
\end{align*}

\end{lemma}

\begin{proof}
(of Theorem \ref{t3.2}). Our present discussion is a tiny generalization of
Proposition 2.7 in \cite{n1}. We define six $(A\times B\times C,(a_{0}%
,b_{0},c_{0}))$-$3$-icons on $M$ as follows:
\begin{align*}
\xi_{123}  &  =\xi_{1}\circledast\xi_{2}\circledast\xi_{3}\\
\xi_{132}  &  =\xi_{1}\circledast(\xi_{2}\widetilde{\circledast}\xi_{3})\\
\xi_{213}  &  =(\xi_{1}\widetilde{\circledast}\xi_{2})\circledast\xi_{3}\\
\xi_{231}  &  =\xi_{1}\widetilde{\circledast}(\xi_{2}\circledast\xi_{3})\\
\xi_{312}  &  =(\xi_{1}\circledast\xi_{2})\widetilde{\circledast}\xi_{3}\\
\xi_{321}  &  =\xi_{1}\widetilde{\circledast}\xi_{2}\widetilde{\circledast}%
\xi_{3}%
\end{align*}
Then it is easy, by dint of Lemma \ref{t3.3}, to see that
\begin{align}
\left\lfloor \xi_{1},\left\lfloor \xi_{2},\xi_{3}\right\rfloor \right\rfloor
&  =(\xi_{123}\underset{1}{\overset{\cdot}{-}}\xi_{132})\overset{\cdot}{-}%
(\xi_{231}\underset{1}{\overset{\cdot}{-}}\xi_{321})\label{J1}\\
\left\lfloor \xi_{2},\left\lfloor \xi_{3},\xi_{1}\right\rfloor \right\rfloor
&  =(\xi_{231}\underset{2}{\overset{\cdot}{-}}\xi_{213})\overset{\cdot}{-}%
(\xi_{312}\underset{2}{\overset{\cdot}{-}}\xi_{132})\label{J2}\\
\left\lfloor \xi_{3},\left\lfloor \xi_{1},\xi_{2}\right\rfloor \right\rfloor
&  =(\xi_{312}\underset{3}{\overset{\cdot}{-}}\xi_{321})\overset{\cdot}{-}%
(\xi_{123}\underset{3}{\overset{\cdot}{-}}\xi_{213})\label{J3}%
\end{align}
Therefore the desired Jacobi identity follows directly from the general Jacobi identity.
\end{proof}

\begin{remark}
In order to see that the right-hand side of (\ref{J1}) is meaningful, we have
to check that all of
\begin{align*}
&  \xi_{123}\underset{1}{\overset{\cdot}{-}}\xi_{132}\\
&  \xi_{231}\underset{1}{\overset{\cdot}{-}}\xi_{321}\\
&  (\xi_{123}\underset{1}{\overset{\cdot}{-}}\xi_{132})\overset{\cdot}{-}%
(\xi_{231}\underset{1}{\overset{\cdot}{-}}\xi_{321})
\end{align*}
are meaningful. Since $\xi_{2}\circledast\xi_{3}\overset{\cdot}{-}\xi
_{2}\widetilde{\circledast}\xi_{3}$ is meaningful by Lemma \ref{t2.3.2},
$\xi_{123}\underset{1}{\overset{\cdot}{-}}\xi_{132}$ is also meaningful and we
have
\[
\xi_{123}\underset{1}{\overset{\cdot}{-}}\xi_{132}=\xi_{1}\circledast(\xi
_{2}\circledast\xi_{3}\overset{\cdot}{-}\xi_{2}\widetilde{\circledast}\xi_{3})
\]
by Lemma \ref{t3.3}. Similarly $\xi_{231}\underset{1}{\overset{\cdot}{-}}%
\xi_{321}$ is meaningful and we have
\[
\xi_{231}\underset{1}{\overset{\cdot}{-}}\xi_{321}=\xi_{1}\widetilde
{\circledast}(\xi_{2}\circledast\xi_{3}\overset{\cdot}{-}\xi_{2}%
\widetilde{\circledast}\xi_{3})
\]
Therefore $(\xi_{123}\underset{1}{\overset{\cdot}{-}}\xi_{132})\overset{\cdot
}{-}(\xi_{231}\underset{1}{\overset{\cdot}{-}}\xi_{321})$ is meaningful by
Lemma \ref{t2.3.2}. Similar considerations apply to (\ref{J2}) and (\ref{J3}).
\end{remark}

\section{The Jacobi Identity for Tangent-Vector-Valued Differential Forms
\label{4}}

Our three distinct but equivalent viewpoints of tangent-vector-valued
differential forms on $M$ are based upon the following exponential laws:
\begin{align*}
\lbrack\lbrack D^{p}  &  \rightarrow M]\rightarrow\lbrack D\rightarrow M]]\\
&  =[[D^{p}\rightarrow M]\times D\rightarrow M]\\
&  =[D\rightarrow\lbrack\lbrack D^{p}\rightarrow M]\rightarrow M]]
\end{align*}
If $p=0$, the above laws degenerate into the corresponding ones in
\S \ref{2.1}.

The first viewpoint, which is highly orthodox, is to regard $[D\rightarrow M]
$\textit{-valued }$p$\textit{-forms on} $M$ as mappings $\omega:[D^{p}%
\rightarrow M]\rightarrow\lbrack D\rightarrow M]$ with $\gamma(0,...,0)=\omega
(\gamma)(0)$ for any $\gamma\in\lbrack D^{p}\rightarrow M]$ and satisfying the
$p$-homogeneity and the alternating \ property in the sense of Definition 1 in
\S 4.1 of Lavendhomme \cite{l1}. By dropping the alternating property, we get
the weakier notion of a $[D\rightarrow M]$\textit{-valued }$p$%
\textit{-semiform on} $M$.

The second viewpoint goes as follows:

\begin{proposition}
\label{t4.8}$[D\rightarrow M]$-valued $p$-forms on $M$ can be identified with
mappings $\omega:[D^{p}\rightarrow M]\times D\rightarrow M$ pursuant to the
following conditions:

\begin{enumerate}
\item $\omega(\gamma,0)=\gamma(0,...,0)$ for any $\gamma\in\lbrack
D^{p}\rightarrow M]$.

\item $\omega(\gamma,\alpha d)=\omega(\alpha\underset{i}{\cdot}\gamma,d)$ for
any $d\in D$, any $\alpha\in\mathbb{R}$, any $\gamma\in\lbrack D^{p}%
\rightarrow M]$ and any natural number $i$ with $1\leq i\leq p$.

\item $\omega(\gamma^{\sigma},d)=\omega(\gamma,\varepsilon_{\sigma}d) $ for
any $d\in D$, any $\gamma\in\lbrack D^{p}\rightarrow M]$ and any $\sigma
\in\mathbb{S}_{p}$.
\end{enumerate}

By dropping the third condition, we get entities corresponding to
$[D\rightarrow M]$-valued $p$-semiforms on $M$.
\end{proposition}

The third viewpoint, which is most radical, goes as follows:

\begin{proposition}
\label{t4.9}$[D\rightarrow M]$-valued $p$-forms on $M$ can be identified with
mappings $\omega:D\rightarrow\lbrack\lbrack D^{p}\rightarrow M]\rightarrow M]$
satisfying the following conditions:

\begin{enumerate}
\item $\omega(0)=\delta_{(0,...,0)}$

\item $\alpha\underset{i}{\cdot}\omega(d)=\omega(\alpha d)$ for any $d\in D$,
any $\alpha\in\mathbb{R}$ and any natural number $i$ with $1\leq i\leq p$.

\item $(\omega(d))^{\sigma}=\omega(\varepsilon_{\sigma}d)$ for any $d\in D$
and any $\sigma\in\mathbb{S}_{p}$.
\end{enumerate}

By dropping the third condition, we get entities corresponding to
$[D\rightarrow M]$-valued $p$-semiforms on $M$.
\end{proposition}

\begin{remark}
By dropping the second and third conditions, we find the notion of
$(D^{p},(0,...,0))$-$1$-icon on $M$.
\end{remark}

The following proposition is simple but very important.

\begin{proposition}
The addition for $[D\rightarrow M]$-\textit{valued }$p$\textit{-semiforms on}
$M$ in the first sense (i.e., using the fiberwise addition of the vector
bundle $[D\rightarrow M]\rightarrow M$) and that in the third sense (i.e., as
the addition of tangent vectors to the microlinear space $[[D^{p}\rightarrow
M]\rightarrow M]$ at $\delta_{(0,...,0)}$) coincide.
\end{proposition}

\begin{proof}
This follows mainly from the exponential law
\begin{align*}
\lbrack\lbrack D^{p}  &  \rightarrow M]\rightarrow\lbrack D(2)\rightarrow
M]]\\
&  =[D(2)\rightarrow\lbrack\lbrack D^{p}\rightarrow M]\rightarrow M]]
\end{align*}
The details can safely be left to the reader.
\end{proof}

Unless stated to the contrary, we will use the terms $[D\rightarrow
M]$-\textit{valued }$p$\textit{-semiforms on} $M$ and $[D\rightarrow
M]$-\textit{valued }$p$\textit{-forms on} $M$ in the third sense.

The following lemma should be obvious.

\begin{lemma}
\label{t4.1}If $\omega_{1}$ is a $[D\rightarrow M]$-valued $p$-semiform on $M
$ and $\omega_{2}$ is a $[D\rightarrow M]$-valued $q$-semiform on $M$, then we
have
\begin{align*}
\alpha\underset{i}{\cdot(}(\omega_{1}\circledast\omega_{2})(d_{1},d_{2}))  &
=(\omega_{1}\circledast\omega_{2})(\alpha d_{1},d_{2})\\
\alpha\underset{i}{\cdot(}(\omega_{1}\widetilde{\circledast}\omega_{2}%
)(d_{1},d_{2}))  &  =(\omega_{1}\widetilde{\circledast}\omega_{2})(\alpha
d_{1},d_{2})
\end{align*}
for any $(d_{1},d_{2})\in D^{2}$ and natural number $i$ with $1\leq i\leq p$,
while we have
\begin{align*}
\alpha\underset{i}{\cdot(}(\omega_{1}\circledast\omega_{2})(d_{1},d_{2}))  &
=(\omega_{1}\circledast\omega_{2})(d_{1},\alpha d_{2})\\
\alpha\underset{i}{\cdot(}(\omega_{1}\widetilde{\circledast}\omega_{2}%
)(d_{1},d_{2}))  &  =(\omega_{1}\widetilde{\circledast}\omega_{2}%
)(d_{1},\alpha d_{2})
\end{align*}
for any $(d_{1},d_{2})\in D^{2}$ and any natural number $i$ with $p+1\leq
i\leq p+q$.
\end{lemma}

\begin{corollary}
\label{t4.2}If $\omega_{1}$ is a $[D\rightarrow M]$-valued $p$-semiform on $M
$ and $\omega_{2}$ is a $[D\rightarrow M]$-valued $q$-semiform on $M$, then
$\left\lfloor \omega_{1},\omega_{2}\right\rfloor $ is a $[D\rightarrow
M]$-valued $(p+q)$-semiform on $M$.
\end{corollary}

\begin{proof}
It suffices to see that
\[
\alpha\underset{i}{\cdot(}\left\lfloor \omega_{1},\omega_{2}\right\rfloor
(d))=\left\lfloor \omega_{1},\omega_{2}\right\rfloor (\alpha d)
\]
for any $d\in D$, any $\alpha\in\mathbb{R}$ and any natural number $i$ with
$1\leq i\leq p+q$, which follows easily from the \ above lemma and Proposition
5 in \S 3.4 of Lavendhomme \cite{l1}.
\end{proof}

Given a $[D\rightarrow M]$-valued $p$-semiform $\omega$\ on $M$ and $\sigma
\in\mathbb{S}_{p}$, we define a $[D\rightarrow M]$-valued $p$-semiform
$\omega^{\sigma}$ to be
\[
\omega^{\sigma}(d)=\omega(d)^{\sigma}%
\]

Now we are ready to state that

\begin{theorem}
\label{t4.3}If $\omega_{1}$ is a $[D\rightarrow M]$-valued $p$-semiform on $M
$ and $\omega_{2}$ is a $[D\rightarrow M]$-valued $q$-semiform on $M$, then we
have
\[
\left\lfloor \omega_{1},\omega_{2}\right\rfloor =-\left\lfloor \omega
_{2},\omega_{1}\right\rfloor ^{\rho}%
\]
where $\rho$ is the permutation mapping the sequence $1,...,q,q+1,...,p+q$ to
the sequence $q+1,...,p+q,1,...,q$.
\end{theorem}

\begin{proof}
This follows simply from Theorem \ref{t3.1}.
\end{proof}

\begin{theorem}
\label{t4.4}If $\omega_{1}$ is a $[D\rightarrow M]$-valued $p$-semiform on $M
$, $\omega_{2}$ is a $[D\rightarrow M]$-valued $q$-semiform on $M$ and
$\omega_{3}$ is a $[D\rightarrow M]$-valued $r$-semiform on $M$, then the
following Jacobi identity holds for the three $[D\rightarrow M]$-valued
$(p+q+r)$-semiforms on $M$:
\[
\left\lfloor \omega_{1},\left\lfloor \omega_{2},\omega_{3}\right\rfloor
\right\rfloor +\left\lfloor \omega_{2},\left\lfloor \omega_{3},\omega
_{1}\right\rfloor \right\rfloor ^{\rho_{1}}+\left\lfloor \omega_{3}%
,\left\lfloor \omega_{1},\omega_{2}\right\rfloor \right\rfloor ^{\rho_{2}}=0
\]
where $\rho_{1}$ and $\rho_{2}$ are the following permutations:
\begin{align*}
\rho_{1}  &  =\left(
\begin{array}
[c]{ccccccccc}%
1 & ... & q & q+1 & ... & q+r & q+r+1 & ... & p+q+r\\
p+1 & ... & p+q & p+q+1 & ... & p+q+r & 1 & ... & p
\end{array}
\right) \\
\rho_{2}  &  =\left(
\begin{array}
[c]{ccccccccc}%
1 & ... & r & r+1 & ... & p+r & p+r+1 & ... & p+q+r\\
p+q+1 & ... & p+q+r & 1 & ... & p & p+1 & ... & p+q
\end{array}
\right)
\end{align*}

\end{theorem}

\begin{proof}
The desired Jacobi identity is a direct consequence of Theorem \ref{t3.2}.
\end{proof}

Now we turn to forms. Given a $[D\rightarrow M]$-valued $p$-semiform $\omega
$\ on $M$ and $\sigma\in\mathbb{S}_{p}$, we define a $[D\rightarrow M]$-valued
$p$-semiform $\omega^{\sigma}$ on $M$ to be
\[
\omega^{\sigma}(d)=\omega(d)^{\sigma}%
\]
Given a $[D\rightarrow M]$-valued $p$-semiform $\omega$\ on $M$, we define a
$[D\rightarrow M]$-valued $p$-semiform $\mathcal{A}\omega$\ on $M$ to be
\[
\mathcal{A}\omega=\sum_{\sigma\in\mathbb{S}_{p}}\varepsilon_{\sigma}%
\omega^{\sigma}%
\]
We write $\mathcal{A}_{p,q}\omega$ for $(1/p!q!)\mathcal{A}\omega$ in case
that $\omega$ is a $[D\rightarrow M]$-valued $(p+q)$-semiform on $M$. We write
$\mathcal{A}_{p,q,r}\omega$ for $(1/p!q!r!)\mathcal{A}\omega$ in case that
$\omega$ is a $[D\rightarrow M]$-valued $(p+q+r)$-semiform on $M$.

Given a $[D\rightarrow M]$-valued $p$-form $\omega_{1}$ on $M$ and a
$[D\rightarrow M]$-valued $q$-form $\omega_{2}$ on $M$, we are going to define
their \textit{Fr\"{o}licher-Nijenhuis bracket} $\left\lceil \omega_{1}%
,\omega_{2}\right\rceil $ to be
\[
\left\lceil \omega_{1},\omega_{2}\right\rceil =\mathcal{A}_{p,q}%
\mathcal{(}\left\lfloor \omega_{1},\omega_{2}\right\rfloor )
\]
which is undoubtedly a $[D\rightarrow M]$-valued $(p+q)$-form on $M$.

\begin{theorem}
\label{t4.5}If $\omega_{1}$ is a $[D\rightarrow M]$-valued $p$-form on $M$ and
$\omega_{2}$ is a $[D\rightarrow M]$-valued $q$-form on $M$, then we have
\[
\left\lceil \omega_{1},\omega_{2}\right\rceil =-(-1)^{pq}\left\lceil
\omega_{2},\omega_{1}\right\rceil
\]

\end{theorem}

\begin{proof}
We have
\begin{align*}
&  \left\lceil \omega_{1},\omega_{2}\right\rceil \\
&  =\mathcal{A}_{p,q}\mathcal{(}\left\lfloor \omega_{1},\omega_{2}%
\right\rfloor )\\
&  =-\mathcal{A}_{p,q}\mathcal{(}\left\lfloor \omega_{2},\omega_{1}%
\right\rfloor ^{\rho})\text{ \ \ [By Theorem \ref{t4.3}]}\\
&  =-\frac{1}{p!q!}\sum_{\tau\in\mathbb{S}_{p+q}}\varepsilon_{\tau}\left(
\left\lfloor \omega_{2},\omega_{1}\right\rfloor ^{\rho}\right)  ^{\tau}\\
&  =-\frac{1}{p!q!}\sum_{\tau\in\mathbb{S}_{p+q}}\varepsilon_{\tau
}\left\lfloor \omega_{2},\omega_{1}\right\rfloor ^{\tau\rho}\\
&  =-\frac{1}{p!q!}\varepsilon_{\rho}\sum_{\tau\in\mathbb{S}_{p+q}}%
\varepsilon_{\tau\rho}\left\lfloor \omega_{2},\omega_{1}\right\rfloor
^{\tau\rho}\\
&  =-\varepsilon_{\rho}\left\lceil \omega_{2},\omega_{1}\right\rceil
\end{align*}
Since $\varepsilon_{\rho}=(-1)^{pq}$, the desired conclusion follows.
\end{proof}

\begin{lemma}
\label{t4.6}If $\omega_{1}$ is a $[D\rightarrow M]$-valued $p$-form on $M$,
$\omega_{2}$ is a $[D\rightarrow M]$-valued $q$-form on $M$ and $\omega_{3} $
is a $[D\rightarrow M]$-valued $r$-form on $M$, then we have
\[
\mathcal{A}_{p,q+r}(\left\lfloor \omega_{1},\mathcal{A}_{q,r}(\left\lfloor
\omega_{2},\omega_{3}\right\rfloor )\right\rfloor )=\mathcal{A}_{p,q,r}%
(\left\lfloor \omega_{1},\left\lfloor \omega_{2},\omega_{3}\right\rfloor
\right\rfloor )
\]

\end{lemma}

\begin{proof}
By the same token as in the familiar associativity of wedge products in
differential forms.
\end{proof}

\begin{theorem}
\label{t4.7}If $\omega_{1}$ is a $[D\rightarrow M]$-valued $p$-form on $M$,
$\omega_{2}$ is a $[D\rightarrow M]$-valued $q$-form on $M$ and $\omega_{3} $
is a $[D\rightarrow M]$-valued $r$-form on $M$, then the following graded
Jacobi identity holds for the three $[D\rightarrow M]$-valued $(p+q+r) $-forms
on $M$:
\[
\left\lceil \omega_{1},\left\lceil \omega_{2},\omega_{3}\right\rceil
\right\rceil +(-1)^{p(q+r)}\left\lceil \omega_{2},\left\lceil \omega
_{3},\omega_{1}\right\rceil \right\rceil +(-1)^{r(p+q)}\left\lceil \omega
_{3},\left\lceil \omega_{1},\omega_{2}\right\rceil \right\rceil =0
\]

\end{theorem}

\begin{proof}
By the same token as in the proof of Theorem \ref{t4.5}. This follows mainly
from Theorem \ref{t4.4} with the help of Lemma \ref{t4.6}\ and the simple fact
that $\varepsilon_{\rho_{1}}=(-1)^{p(q+r)}$ and $\varepsilon_{\rho_{2}%
}=(-1)^{r(p+q)}$. The details can safely be left to the reader.
\end{proof}

\section{The Jacobi Identity for Schwartz Distributions\label{5}}

By a \textit{distribution with compact support on} $M$ we mean a mapping
$u:[M\rightarrow\mathbb{R]\rightarrow R}$ with the property that
\[
u(\alpha f)=\alpha u(f)
\]
for any $f\in\lbrack M\rightarrow\mathbb{R]}$ and any $\alpha\in\mathbb{R}$.
Let us suppose that we are given $x\in M$. By a \textit{Dirac }$x$%
\textit{-flow on} $M$, we mean a $(M,x)$-$1$-icon $\xi$ on $\mathbb{R}$ with
the property that $\xi(d)$ is a distribution with compact support on $M $ for
any $d\in D$.

\begin{lemma}
\label{t5.1}If $u$ is a distribution with compact support on $M$ and $v$ is a
distribution with compact support on $N$, then $u\ast v$ as well as
$u\widetilde{\ast}v$ is a distribution with compact support on $M\times N$.
\end{lemma}

\begin{proof}
We note that, given $\alpha\in\mathbb{R}$ and $h\in\lbrack M\times
N\rightarrow\mathbb{R]}$, we have
\begin{align*}
&  (u\ast v)(\alpha h)\\
&  =u(\lambda x\in M.v(\lambda y\in N.\alpha h(x,y)))\\
&  =u(\lambda x\in M.v(\alpha(\lambda y\in N.h(x,y))))\\
&  =\alpha(u(\lambda x\in M.v(\lambda y\in N.h(x,y))))\\
&  =\alpha(u\ast v)(h)
\end{align*}
so that $u\ast v$ is a distribution with compact support on $M\times N$.
Similarly for $u\widetilde{\ast}v$.
\end{proof}

\begin{proposition}
If $\xi_{1}$ is a Dirac $x$-flow on $M$ and $\xi_{2}$is a Dirac $y$-flow on
$N$ , then $\left\lfloor \xi_{1},\xi_{2}\right\rfloor $ is a Dirac
$(x,y)$-flow on $M\times N$.
\end{proposition}

\begin{proof}
It suffices to note that the space of distributions with compact support on
$M$ forms a microlinear space, from which the desired result follows from the
above Lemma.
\end{proof}

\begin{theorem}
If $\xi_{1}$ is a Dirac $x_{1}$-flow on $M_{1}$, $\xi_{2}$is a Dirac $x_{2}
$-flow on $M_{2}$ and $\xi_{3}$ is a Dirac $x_{3}$-flow on $M_{3}$, then we
have
\[
\left\lfloor \xi_{1},\left\lfloor \xi_{2},\xi_{3}\right\rfloor \right\rfloor
+\left\lfloor \xi_{2},\left\lfloor \xi_{3},\xi_{1}\right\rfloor \right\rfloor
+\left\lfloor \xi_{3},\left\lfloor \xi_{1},\xi_{2}\right\rfloor \right\rfloor
=0
\]

\end{theorem}

\begin{proof}
This is a direct consequence of Theorem \ref{t3.2}.
\end{proof}


\begin{thebibliography}{99}                                                                                               %
\bibitem {ama}Amalio, Roberto M. and Curien, Pierre-Louis:Domains and
Lambda-Calculi, Cambridge University Press, Cambridge, 1998.

\bibitem {bar}Barr, Michael and Wells, Charles:Category Theory for Computer
Science, Prentice Hall, London, 1990.

\bibitem {fn}Fr\"{o}licher, A. and Nijenhuis, A.:Theory of vector-valued
differential forms, Part I, Indagationes Math., \textbf{18} (1956), 338-359.

\bibitem {k0}Kakita, T.:An Introduction to Schwartz Distributions (2nd
edition), in Japanese, Nihonhyoronsha, Tokyo, 1999.

\bibitem {kl}Kock, A. and Lavendhomme, R.:Strong infinitesimal linearity, with
applications to strong difference and affine connections, Cahiers de Topologie
et Geom\'{e}trie Differ\'{e}ntielle, \textbf{25} (1984), 311-324.

\bibitem {k2}Kock, A.:Synthetic Differential Geometry (2nd edition), Cambridge
University Press, Cambridge, 2006.

\bibitem {l1}Lavendhomme, R.: Basic Concepts of Synthetic Differential
Geometry, Kluwer, Dordrecht, 1996.

\bibitem {m0}Michor, Peter W.:Remarks on the Fr\"{o}licher-Nijenhuis bracket,
Proceedings of the Conference on Differential Geometry and its Applications,
Brno 1986, D. Reidel, 1987, pp.197-220.

\bibitem {m1}Michor, Peter W.:Topics in Differential Geometry, American
Mathematical Society, Providence, Rhode Islands, 2008.

\bibitem {m2}Minguez, M.C.:Wedge products of forms in synthetic differential
geometry, Cahiers de Topologie et Geom\'{e}trie Differ\'{e}ntielle,
\textbf{29} (1988), 59-66.

\bibitem {m3}Minguez, M.C.:Some combinatorial calculus on Lie derivatives,
Cahiers de Topologie et Geom\'{e}trie Differ\'{e}ntielle, \textbf{29} (1988), 241-247.

\bibitem {n0}Nijenhuis, A.:Jacobi type identities for bilinear differential
concomitants of certain tensor fields I, Indagationes Math., \textbf{17}
(1955), 390-403.

\bibitem {n1}Nishimura, H.:Theory of microcubes, International Journal of
Theoretical Physics, \textbf{36} (1997), 1099-1131.

\bibitem {n2}Nishimura, H.:General Jacobi identity revisited, International
Journal of Theoretical Physics, \textbf{38} (1999), 2163-2174.

\bibitem {n3}Nishimura, H. and Osoekawa, K.:General Jacobi identity revisited
again, International Journal of Theoretical Physics, \textbf{46} (2007), 2843-2862.

\bibitem {n4}Nishimura, H.:The Fr\"{o}licher and Nijenhuis calculus in
synthetic differential geometry, ArXiv 0810.5492.
\end{thebibliography}
\end{document}